\documentclass{amsart}
\usepackage{graphicx}    
\usepackage{amssymb, mathrsfs, amscd}
\usepackage{Mytheorems}
\author{Alexander Sotirov}
\address{Stanford University} 
\email {sotirov@math.stanford.edu}

\begin{document}     
\title{The Boundary Structure of Zero-Temperature Driven Hard Spheres}

\begin{abstract}
We study the fundamental problem of two gas species whose molecules collide as hard spheres in the presence of a flat boundary and with dependence on only one space dimension. More specifically the steady linear problem considered is the one arising when the second gas dominates as a flow moving towards the boundary with constant microscopic velocity (and hence zero temperature). The boundary condition adopted consists of prescribing the outgoing velocity distribution of the first gas at the boundary. It is discovered that the presence of the boundary under general assumptions on the outgoing distribution ensures the convergence of a series of path integrals resulting in a convenient representation for the distribution of the velocities of the molecules returning at the boundary.
 
\end{abstract} 

\maketitle
 
\begin{section}{Introduction and Statement of the Result}

We begin by noting that the Boltzmann equation with boundary conditions has been previously considered in a number of works, for example the reader may consult the work by Bardos, Caflisch, and Nicolaenko \cite{bardos}, the zero-temperature paper by Caflisch \cite{caflisch} where behavior in the interior is studied in a perturbative setting with regard to applications to strong shock waves, as well as the works \cite{golse} and \cite{sone1} among many others. Our approach differs substantially from the previous works since we investigate the structure of the solution in terms of contributions from the individual particle paths and we look for an exact description of the solution at the boundary in the velocity variable. Excellent reviews of kinetic theory in general are \cite{cerc} and \cite{villani} as well as \cite{sone}.

In our problem we consider molecules emitted from a flat boundary which interact as hard spheres with the molecules of a second gas which have zero temperature and penetrate into the boundary. We will be interested in the distribution of the velocities of the emitted molecules at the time of their return to the boundary. The interactions between the emitted molecules themselves are ignored and so a linear problem is studied. Accounting for the nonlinear interactions would be an important continuation of this work. We must also point out that the zero temperature distribution of the background is of course an idealization. However it is reasonable to expect that many of the features studied here will be still present in the case of the more natural Gaussian distribution of the background (see the concluding remarks).

The equation at hand is time independent, valid for $x\geq 0$, and takes the form
\begin{align*} \xi^1 \partial_x f(\xi,x) &=  Q(f,\rho \delta(\xi - {\bf c}))\\
f(\xi, 0)&=f^{+}(\xi) \text { for } \xi^1>0
\end {align*}
where $f(\xi,x)$ with $\xi=(\xi^1,\xi^2, \xi^3)$ is the unknown velocity distribution and $\rho\delta(\xi - {\bf c}) $ with ${\bf c}=(-c,0,0)$ for $c>0$  is the background distribution where $\rho$ is the number density of the background.  So the background molecules enter the boundary and their distribution is homogeneous in the interior. Here $Q$ is the Boltzmann collision operator. To simplify the presentation let us assume that both particle species have the same molecular mass, normalized to 1. The case when the masses differ can be treated similarly. We will focus on determining the value $f^{-}(\xi)=f(\xi,0)$ for $\xi^1<0$.  
After writing out $Q$ as the difference between the gain and loss terms due to collisions we have
\begin {align*}
\xi^1 \partial_x f(\xi,x)  &= \rho \sigma^2 \int_{\mathbb{R}^3}\int_{S^+} f(\xi-((\xi-\xi_*)\cdot n) n)\delta(\xi_*+((\xi-\xi_*)\cdot n) n - {\bf c})|(\xi-\xi_*)\cdot n| dn d\xi_* \\
&- \rho\sigma^2 \int_{\mathbb{R}^3} \int_{S^+} f(\xi)\delta(\xi_*-{\bf c})|(\xi_*-\xi)\cdot n| dn d\xi_* =\\
&=\rho\sigma^2 Kf - \rho\sigma^2\pi|\xi-{\bf c}|f.
\end{align*}
where $n$ is the unit vector at which the collision occurs and $S^{+}$ is the half sphere $(\xi-\xi_*)\cdot n >0$, $\sigma$ is the sphere radius. For our purpose we will take $\rho\sigma^2 =1$. This is actually no restriction as explained in the remark on scaling after the statement of the theorem. 
We will give a more convenient expression for $Kf(\xi_0)$, but first to elucidate it we will determine the region containing particles $\xi_1$ which can influence the region near $\xi_0$ after experiencing a collision at the direction of $n$. In other words given $\xi_0=(\xi_0^1,\xi_0^2, \xi_0^3)$ we need to determine the set $\{n, \xi_1 =(\xi_1^1,\xi_1^2, \xi_1^3 )\}$ such that 
$$\xi_0= \xi_1 - ((\xi_1- {\bf c})\cdot n) n.$$  
In the coordinate system $n, n^\perp$ let $P_{n}$ and $P_{n}^{\perp}$ denote the corresponding projections. Since the effect of a collision is to exchange the momentum along $n$, the new momentum along $n$ of the particle hit by a background particle is always $P_{n}({\bf c})$, i.e. we see that $P_{n}({\bf c})=P_n(\xi_0)$. This is satisfied iff $n \perp \xi_0 - {\bf c}$ which determines $n$. But then we must also have $P^{\perp}_{n}(\xi_1)=P^{\perp}_{n}(\xi_0)$. Therefore the particles with velocity $\xi_1$ that can influence the velocity $\xi_0$ lie along the plane $L_{{\bf c} \xi_0}$ passing through $\xi_0$ and perpendicular to $\xi_0-{\bf c}$ (see Fig. 1). In fact it can be easily shown that the interior of the sphere $C_{\xi_0}$ on the figure is the region influenced by $\xi_0$ after an arbitrary number of collisions and the exterior is the region that influences $\xi_0$. This answers our question.
\begin {figure}
\begin{center}
\includegraphics[width = 8 cm]{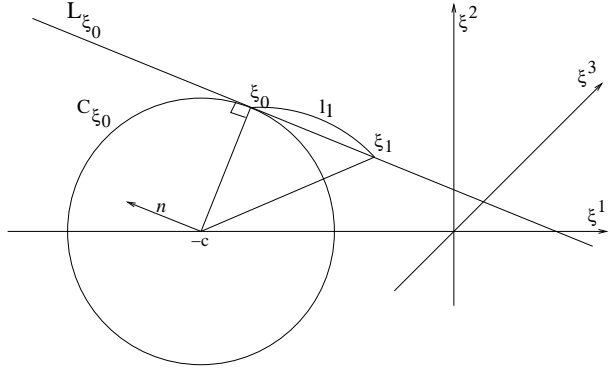}  
\end{center}
\caption{The plane of influence}
\end{figure}
So clearly $K$ will be given by an integral along this plane but we need to determine the correct weight. Although a direct calculation is of course possible, for brevity to obtain the exact formula we will make use of Torsten Carleman's well-known representation (\cite{carl1}, \cite{carl2}) for the gain term, namely:
$$ Q^+(f,g)(\xi) = \int_{\mathbb{R}^3} \int_{S^+} f(\xi-((\xi-\xi_*)\cdot n) n)g(\xi_*+((\xi-\xi_*)\cdot n) n)|(\xi-\xi_*)\cdot n| dn d\xi_* = $$
$$= \int_{\mathbb{R}^3} \frac{g(\xi_1)}{|\xi-\xi_1|} \int_{L_{\xi_1,\xi}}f(\xi_2)d\xi_2 d\xi_1,$$
where $L_{\xi,\xi_1}$ denotes the plane perpendicular to $\xi-\xi_1$ and passing through $\xi$. For completeness the proof of this representation is briefly indicated in the Appendix. Specializing the above to the case which interests us, namely $g(\xi)=\delta(\xi -{\bf c})$ we obtain that
\begin{equation} \label {kform} Kf(\xi_0) = \frac{1}{|\xi_0-{\bf c}|}\int_{L_{{\bf c}, \xi_0}} f(\xi_1) d\xi_1. \end{equation}
In what follows we will omit ${\bf c}$ from the subscript of the plane since ${\bf c}$ is fixed.

We are going to prove the following 

\begin{theorem} \label{theo} Let there be given an outgoing (from the boundary) distribution $f^{+}(\xi)$ satisfying $0 \leq f^{+} < M$ which is supported in the compact region $\xi^{1}>0$ and $|\xi-{\bf c}|<R$. Then for the incoming (to the boundary) distribution denoted by $f^{-}(\xi_0)$ and defined for $\xi_0^1 \leq 0$ we have the following:

\noindent i) An explicit series with positive terms with factorial convergence exists such that
$$f^-(\xi_0) = \sum_{n=1}^{\infty} f_n^-(\xi_0)$$
where the $n$-th term represents the contribution of particles returning after $n$ collisions, and specifically the rate of convergence is such that if $\mathcal{R}_n$ denotes the remainder term then the mass flux due to it satisfies
$$\int_{\mathbb{R}^3_{\xi_0},\xi_0^1<0} |\xi_0^1|\mathcal{R}_n (\xi_0) d\xi_0 < \frac{MR^5}{c}ne^{-n+1}.$$
whenever it is satisfied that 
$$n > 30(\log\frac{R}{c}+10).$$
We see that the number of terms describing the solution grows linearly in $\log\frac{R}{c}$.  

\noindent ii) If we use the notation $r_0=|\xi_0-{\bf c}|$ the solution $f^-_n(\xi_0)$ is identically $0$ for $r_0>R$ and it possesses a singularity at $r_0=0$. The specific nature of the singularity is such that for $r_0 < c/2$ and any $\delta >0$ we have 
$$f^{-}(\xi_0) \leq  C_{\delta} (c,R)\frac{M R^3}{c  |\xi_0^1|} \frac{1}{r_0^{1+\delta}}$$
were $C_{\delta}$ is an explicit constant whose dependence on $c$ and $R$ is indicated in the proof.
\end{theorem}

We remark that although the above bound already indicates fast convergence, a more careful analysis may improve notably the constants given. But for us it is most important that these constants are explicit and that in our proof we establish the mechanism leading to this convergence. The estimates indicated above present integrable functions and the first terms of the series mentioned above contain essentially all the mass flux of the solution. The mass flux is the natural measure of the distribution at the boundary (as is made precise for example in \cite{cessenat}). The results of this theorem are in contrast to other situations in kinetic theory where such expansions cannot be expected to converge with any reasonable rate. In these other situations the lack of fast convergence is due to the fact that (fluid) components of the solution corresponding to the kernel elements of the collision operator (and representing thermodynamic equilibrium) propagate in the interior of the domain of interest (in our case $K-\nu$ has only one kernel element which is the Dirac delta and our linear equation only has conservation of mass). In our problem the characteristics do not propagate inside the domain so these fluid components are irrelevant and we are studying a purely kinetic phenomenon. To study fluid components propagating in the interior (at least in the case of Gaussian background where $L^2$ space methods are applicable) understanding of the spectrum and the use of perturbation theory are most natural (see for instance \cite{LiuYu}). Although we do not investigate this here, it is reasonable to expect that a description similar to ours will also apply to the kinetic part of the solution in the case of the so called linearized collision operator (when a perturbation of the background is considered so energy and momentum are also conserved) but in such a case depending on the speed of the background and the conditions at infinity a fluid component of the solution may also appear as a remainder term in the above series.

The behavior $r_0 \rightarrow 0$ presented by the above formula is essentially optimal, in fact except for the $\delta$-correction it is already displayed by the first term of the series. In the limit $c \rightarrow 0$ we have $C_{\delta}/c \rightarrow \infty$. This reflects the fact that when the background velocity vanishes the linearized problem does not posses a finite steady solution with zero flux at infinity. The equilibrium in the case $c \approx 0$ is due to the particles emitted from the wall interacting with themselves which is a nonlinear effect, or alternatively a linearization is necessary around certain density which does not vanish away from the wall.

We hope that the understanding gained from this problem may be useful in studying other boundary conditions, as for example the idealized situation of specular reflection $f(\xi^1,\xi^2,\xi^3) = f(-\xi^1,\xi^2,\xi^3)$ at $x=0$. In this problem it would be quite interesting to know the momentum flux $\int_{\mathbb{R}^3(\xi)}(\xi^1)^2f(\xi)d\xi$ at $x=0$ which gives the force exerted on the boundary by the background through the interaction with the particles trapped near the boundary.

Now we briefly discuss the scaling in this problem. There is no macroscopic reference length in the $x$ variable, so a Knudsen number cannot be specified. Furthermore since we are taking $x=0$ as the region of attention we cannot expect our solution to be influenced by the dimensional factor $\rho\sigma^2$. This is in fact exactly so: the above theorem would hold unchanged even if there was such a factor in the right hand side of the original equation. The reason is that we are not interested in how long or how far a particle had to travel before return, we only care for the velocity at return. A proof of this is noted after equation \eqref{scale}. 
\end{section}

\begin{section}{The Particle Path Representation and Proof of the Result}

To prove the above result we will consider the contributions of the different paths a particle coming out of the boundary can travel until it returns to the wall as an incoming particle. More specifically we will distinguish these paths according to the number of collisions experienced before return. An easy approach to this representation comes from the Laplace transform formalism. Let $\mathcal{L}f(\xi, z)$ denote the Laplace transform in the $x$ variable. We have:
$$ \xi^1 \partial_x f(\xi, x) = Kf - \nu f \Rightarrow$$
$$ \xi^1 [z \mathcal{L}(f) - f(\xi, 0) ] = \mathcal{L}(Kf) - \nu \mathcal{L}(f) \Rightarrow$$
$$ \mathcal{L}(f) = [z\xi^1 - K + \nu]^{-1} \xi^1 f(\xi,0) \Rightarrow$$
$$ f(\xi, x) = \mathcal{L}^{-1}\left([z\xi^1-K+\nu]^{-1}\xi^1f(\xi, 0) \right) \Rightarrow$$
$$f(\xi, x) = \frac{1}{2\pi i}\int_{-i\infty}^{i\infty} e^{zx}[z\xi^1-K+\nu]^{-1}\xi^1 f(\xi, 0) dz,$$
or if we switch to the more convenient Fourier notation
$$f(\xi, x) = \frac{1}{2\pi}\int_{-\infty}^{\infty} e^{ikx}[ik\xi^1-K+\nu]^{-1}\xi^1 f(\xi, 0) dk$$
Now if we write 
$$[ik\xi^1-K+\nu]^{-1} = [ik\xi^1+\nu]^{-1}\left[I-K[ik\xi^1+\nu]^{-1}\right]^{-1}$$
and introduce the notation $\Phi(k,\xi)=K[ik\xi^1+\nu]^{-1}$ we can write the expansion 
$$[I-\Phi(k,\xi)]^{-1} = I + \Phi(k,\xi) + \Phi(k,\xi)^2 + \ldots$$
and we obtain a series representation for the solution. It is the convergence of the series resulting from the above representation that we shall establish, provided the outgoing distribution satisfies the requirements in the theorem. The zeroth term in the series represents particles experiencing no collision at all after leaving the wall, so it will be ignored. The $n-th$ term's contribution to the incoming distribution is 
\begin{align} \label{scale}
&f^{-}_n(\xi_0) = \frac{1}{2\pi}{{\int_{-\infty}^{\infty} \frac{ e^{ikx} } { ik\xi_0^1+\pi|\xi_0-{\bf c}|}  \Phi^n(k,\xi)\xi^{1}f^{+}(\xi) dk} \hspace {2mm}\vline}_{x=0} = \notag\\
&\frac{1}{2\pi} \int_{-\infty}^{\infty} \frac{1}{ik\xi_0^1+\pi|\xi_0-{\bf c}|}K  \frac{1}{ik\xi_1^1+\pi|\xi_1-{\bf c}|}  \\ 
&K  \frac{1}{ik\xi_2^1+\pi|\xi_2-{\bf c}|} \ldots K \frac{1}{ik\xi_n^1+\pi|\xi_n-{\bf c}|} |\xi_n^1|f^{+}(\xi_n) dk. \notag
\end {align}
and we will construct $f^{-}$ as 
$$f^{-} = \sum_{n=1}^{\infty} f^{-}_n.$$
It is easy to see that this decomposition corresponds to the number of collisions experienced by a particle, since by the convolution theorem for the Fourier transform 
\begin{multline*}
\mathcal{F}^{-1}\left(\frac{1}{\prod_{i=0}^n(ik+\frac{\pi|\xi_i-{\bf c}|}{\xi_i^1})}\right) = \\
 \mathcal{F}^{-1}\left(\frac{1}{ik+\frac{\pi|\xi_0-{\bf c}|}{\xi_0^1}}\right)\star \mathcal{F}^{-1}\left(\frac{1}{ik+\frac{\pi|\xi_1-{\bf c}|}{\xi_1^1}}\right)\star \ldots \star \mathcal{F}^{-1}\left(\frac{1}{ik+\frac{\pi|\xi_n-{\bf c}|}{\xi_n^1}}\right)=\\
=\left(e^{-\frac{\pi|\xi_0-{\bf c}|}{\xi_0^1}x }\chi_{\{\text{sgn}(\xi_0^1)x>0\}}\right)\star \left(e^{-\frac{\pi|\xi_1-{\bf c}|}{\xi_1^1}x} \chi_{\{\text{sgn}(\xi_1^1)x>0\}}\right)\star \ldots \star \left(e^{-\frac{\pi|\xi_n-{\bf c}|}{\xi_n^1}x} \chi_{\{\text{sgn}(\xi_n^1)x>0\}}\right) 
\end{multline*}
where the $\chi$ are just indicator functions. This expression is a convolution of exponential distributions (up to the normalizing constants) where the exponent is $1/l$ where $l$ is the expected distance a particle will travel before it experiences a collision ($1/l$ is the ratio of the collision frequency to the velocity). Our formula \eqref{scale} differs from this by the additional information due to $K$ which specifies what new velocity the particle is likely to acquire after the collision is experienced. With this it becomes clear that iterating $n$ steps represents a path consisting of $n$ collisions.

To address the scaling question mentioned in the introduction, note that if we had included $\rho\sigma^2$ then in \eqref{scale} after changing $k'=k/(\rho\sigma^2)$, $dk = \rho\sigma^2 dk'$ this constant will disappear at $x=0$ as desired. At $x>0$ it of course doesn't disappear but will appear in an exponential factor of spatial decay.
 
Now fix a set of values $\xi_0, \xi_1, \dots, \xi_n$ in the above expression. Let it be the case that $\xi_{for_i}^1>0$ and $\xi_{ret_j}^1<0$ where $\{for_i\}$ and $\{ret_j\}$ are sets of indices which partition $\{0,1,\dots,n\}$ with $i= 1,2,\dots,l$, and $j=1,2,\dots,m$ where $l+m=n+1$. The first correspond to particles moving forwards (away from the wall) and the second to returning particles (moving backwards). Throughout its $n$-collision trip a particle could experience an arbitrary sequence of forward and backward movements. We will perform the $dk$ integration with these fixed velocities. For this it is necessary to compute the following integral (distinguish $\sqrt{-1}$ from the index $i$):
$$ \int_{-\infty}^{\infty} \frac{dk}{\prod_{i=1}^{l}(ik|\xi_{for_i}^1|+\pi|\xi_{for_i}-{\bf c}|) \prod_{j=1}^{m}(-ik|\xi_{ret_j}^1|+\pi|\xi_{ret_j}-{\bf c}|)}.$$ 
Introducing the notation
$$a_i = \frac{\pi|\xi_{for_i}-{\bf c}|}{|\xi_{for_i}^{1}|}, \text {   } b_j = \frac{\pi|\xi_{ret_j}-{\bf c}|}{|\xi_{ret_j}^{1}|} $$
The above integral can be written as
$$ \frac{\prod_{i=1}^l a_i \prod_{j=1}^m b_j} {\pi^{n+1}\prod_{i=0}^n |\xi_i-{\bf c}|} \int_{-\infty}^{\infty} \frac{dk}{\prod_{i=1}^{l}(ik+a_i) \prod_{j=1}^{m}(b_j-ik)}.$$ 
It is possible to evaluate the above integral explicitly. In fact its evaluation is the following Lemma

\begin {lemma} 
\noindent i) If for $a_i>0, b_j>0$ we use the notation
$$ J^{l,m} (a_1,\dots, a_l, b_1, \dots b_m)= \int_{-\infty}^{\infty} \frac{\prod_{i=1}^la_i\prod_{j=1}^m b_j }{\prod_{i=1}^{l}(ik+a_i) \prod_{j=1}^{m}(b_j-ik)} dk$$
then the following recursive relation holds
$$ J^{l,m}= \frac{a_l}{a_l+b_m} J^{l-1,m} + \frac{b_m}{a_l+b_m}J^{l,m-1}.$$

\noindent ii)  Suppose that we have $k$ indices $j_1, j_2, \dots, j_k$ such that $b_{j_s} < B$ for $s = 1\dots k$ and we also have that $a_i>A$ for all $i$. Then we have
$$J^{l,m} \leq \frac{1}{2\pi} B^k \left(1+\frac{1}{A}\right)^{l+m}. $$

\end {lemma}
\begin {proof}
\noindent i) Notice that
$$ a_lJ^{l-1,m}+b_mJ^{l,m-1} = \int_{-\infty}^{\infty} \frac{(ik+a_l + b_m -ik)\prod_{i=1}^la_i\prod_{j=1}^mb_j}{\prod_{i=1}^{l}(ik+a_i) \prod_{j=1}^{m}(b_j-ik)} dk  = (a_l+b_m) J^{l,m}$$
$$\Rightarrow  J^{l,m} =  \frac{a_l}{a_l+b_m} J^{l-1,m} + \frac{b_m}{a_l+b_m}J^{l,m-1}$$
as desired.

\noindent ii) Notice that $J^{l,m}(a_1, \dots, a_l, b_1, \dots, b_m)$ is a symmetric function in each set of variables. So without loss of generality we can assume that $b_1 \leq B$. 

The recurrence relation from i) tells us we can construct $J^{l,m}$ by constructing a tree with a root at $J^{l,m}$ and each node having either one or two children. All the leaves of the tree contain either the expression $\frac{a_2}{a_2+b_1}J^{1,1}$ or the expression $\frac{b_2}{a_1+b_2}J^{1,1}$. This is because by elementary contour integration we have $J^{0,s}=J^{s,0} =0$ for all $s > 1$. Also notice that by contour integration
$$J^{1,1}=\int_{-\infty}^{\infty} \frac{a_1b_1}{(ik+a_1)(-ik+b_1)} dk = \frac{1}{2\pi}\frac{a_1b_1}{a_1+b_1} = \frac{1}{2\pi}\frac{b_1}{1+\frac{b_1}{a_1}} \leq \frac{1}{2\pi}B$$
We will prove the claim by induction on $m,l,k$. By the above it holds for $l=1, m=1, k =1$.

\noindent \textit {Case 1:} $b_m$ not in $\{ b_{j_1}, \dots, b_{j_k} \}$. Then $J^{l-1, m}$ and $J^{l,m-1}$ both still contain $k$ arguments $b_j$ with $b_j<B$, so the induction hypothesis applies to them with $k$ unchanged, i.e.
$$ J^{l,m} = \frac{a_l}{a_l+b_m} J^{l-1,m} + \frac{b_m}{a_l+b_m}J^{l,m-1} \leq \max \left\{J^{l-1, m}, J^{l,m-1}\right\},$$
and the result follows by induction.

\noindent \textit{Case 2:} $b_m \in \{ b_{j_1}, \dots, b_{j_k} \}$. Then $J^{l,m-1}$ has at least $k-1$ arguments $b_j$ with $b_j<B$ and $J^{l-1,m}$ still has $k$ arguments $b_j$ with $b_j<B$. So applying the induction hypothesis
$$ J^{l,m} \leq 
\frac{1}{2\pi} \frac{b_m}{a_l+b_m} B^{k-1}\left(1+\frac{1}{A}\right)^{m+l-1}+
\frac{1}{2\pi} \frac{a_l}{a_l+b_m} B^k \left(1+\frac{1}{A}\right)^{m+l-1} \leq$$
$$\leq \frac{1}{2\pi}  \frac{B^k}{A}\left(1+\frac{1}{A}\right)^{l+m-1}+ \frac{1}{2\pi}B^k\left(1+\frac{1}{A}\right)^{l+m-1} = 
\frac{1}{2\pi}  B^k\left(1+\frac{1}{A}\right)^{l+m}.$$
as desired.
\end  {proof}
\begin{remark} The above estimate is not optimal, but it will be by far sufficient for our purpose. Notice that the expressions $J^{l,m}$ are trivially all less than $J^{1,1}$. In particular the growth due to the $(1+\frac{1}{A})$ term in the above is not optimal. 
\end {remark}
Let $\alpha_i$ be the angle between $\xi_{i-1}$ and $\xi_i$ with respect to a pole centered at ${\bf c}$ and let $r_i=|{\bf c}-\xi_i|$. Let $l_i, \phi_i$ be the polar coordinates in the plane $L_{\xi_{i-1}}$ centered at $\xi_{i-1}$. This notation is illustrated on Fig. 2. With this we can write \eqref{kform} in the more convenient form
$$ Kf(\xi_{i-1}) = \frac{1}{r_{i-1}} \int_0^{2\pi} \int_0^{\infty} l_i f(l_i, \phi_i) dl_i d\phi_i.$$
Eventually we will pick the angles $\alpha_i$ as variables of integration. So we consider the vectors
$$\xi_i(\xi_0, \alpha_1, \phi_1,  \alpha_2, \phi_2 \ldots \alpha_n, \phi_n).$$ 
From the definition we must restrict $\alpha_i \in [0, \pi/2)$. We have that
$$r_i=\frac{r_{i-1}}{\cos \alpha_i} \text { and so } r_i=\frac{r_0}{\cos \alpha_1 \cos \alpha_2 \ldots \cos \alpha_i}$$
Notice that the angles $\phi_1,\phi_2, \dots \phi_i$ are not needed to express $r_i$. This is important in what follows.
\begin{figure}
\begin{center}
\includegraphics[height = 6 cm]{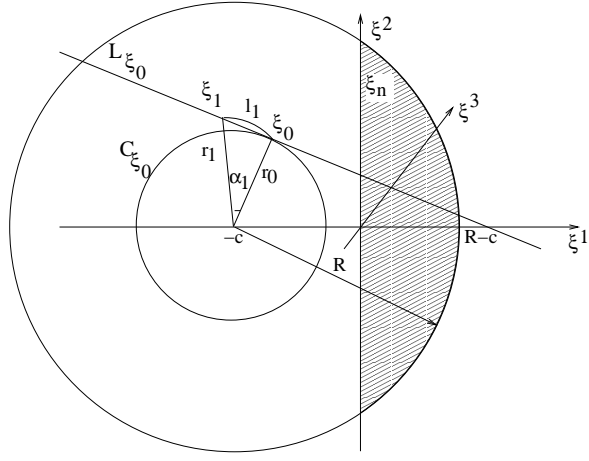} \\
\end{center}
\caption{The support of the initial data and notation}
\end{figure}
Due to the assumption on the support of $f^{+}$ (the shaded area on Fig. 2) we require that $ c\leq r_n\leq R$ which translates to
\begin {equation} \label{mainineq} c \leq r_n=\frac{r_0}{\cos \alpha_1 \cos \alpha_2 \ldots \cos \alpha_n} \leq R. \end {equation}
This restriction is clearly not optimal, since it includes an area larger than the support of $f^+$. This however will only affect the final estimate by a fixed constant related to the ratio of the volumes of the parts of the shell $c<r_n<R$ on the two sides of $\xi^1=0$. In any case this does not affect the convergence structure that we study. 

Let $W$ be the set of $\alpha_i$ in $[0,\pi/2)^n$ which satisfy the above condition \eqref{mainineq}. Also let $U$ denote $[0,2\pi]^n$ which is the domain of integration in the $\phi_i$. We also have
$$l_i = r_{i-1} \tan \alpha_i \text { and } dl_i = \frac{r_{i-1}}{\cos^2 \alpha_i} d\alpha_i.$$
With this notation, the definition of $K$ from \eqref{kform} and using the Lemma we can write:
\begin {align*}  
&f^{-}_n(\xi_0)= \frac{1}{2\pi}\int_{-\infty}^{\infty}\int_U\int_W 
\frac{1}{ik\xi_0^1+\pi|\xi_0-{\bf c}|}\frac{1}{|\xi_0-{\bf c}|}l_1(\xi_0, \alpha_1) \\
&\frac{1}{ik\xi_1^1+\pi|\xi_1(\xi_0,\alpha_1, \phi_1)-{\bf c}|}\frac{1}{|\xi_1(\xi_0,\alpha_1,\phi_2)-{\bf c}|} l_2(\xi_0, \alpha_1,\phi_1 \alpha_2, \phi_2)\\
&\ldots \frac{1}{ik\xi_n^1+\pi|\xi_n-{\bf c}|}\frac{1}{|\xi_n-{\bf c}|} l_n(\xi_0, \alpha_1, \phi_1,  \dots, \alpha_n, \phi_n )
 |\xi_n^1(\xi_0, \alpha_1, \phi_1 \dots, \alpha_n, \phi_n)|\\
&f^{+}(\xi_n)
 \frac{dl_n}{d\alpha_n} \frac{dl_{n-1}}{d \alpha_{n-1}} \ldots \frac{dl_2}{d\alpha_2}\frac{dl_1}{d\alpha_1} d\alpha_n, \ldots d\alpha_1 d\phi_1\dots d\phi_n dk = \\
&= \frac{1}{2\pi^{n+2}}\int_U \int_W \left(\frac {1}{r_0}\frac{\sin\alpha_1}{\cos\alpha_1}r_0\frac{1}{r_1}\frac{\sin\alpha_2}{\cos \alpha_2} r_1 \ldots \frac{1}{r_{n-1}}\frac{\sin \alpha_n}{\cos\alpha_n} r_{n-1} \right) \\
&\left(r_0\frac{1}{\cos^2\alpha_1}r_1\frac{1}{\cos^2\alpha_2}\ldots r_{n-1}\frac{1}{\cos^2\alpha_n} \right) \frac{1}{r_0r_1\dots r_n} \\
&J(a_1,\dots,a_l, b_1, \dots, b_m) |\xi^1_n|f^{+}(\xi_n)  d\alpha_n \ldots d\alpha_2 d\alpha_1 d \phi_n, \ldots d \phi_1 =\\
&=\frac{1}{2\pi^{n+2}}\int_U \int_W \frac{\prod_{i=1}^{n}\sin \alpha_i} {(\prod_{i=1}^{n} \cos \alpha_i)^3} J(a_1,\dots,a_l, b_1, \dots, b_m)\frac{|\xi_n^1|}{r_n}f^{+}(\xi_n)  d\alpha_n \ldots d\alpha_2 d\alpha_1 d\phi_n \ldots d\phi_1.
\end{align*}
where $J(a_1,\dots,a_l, b_1, \dots, b_m)$ denotes the integral from Lemma 2.1. The numbers $l$ and $m$ here are a function of the specific choice $\xi_0, \xi_1, \dots \xi_n$. Now using \eqref {mainineq} we arrive at
\begin{align} \label{fform} 
f^{-}_n(\xi_0) < \frac{1}{2\pi^{n+2}} \int_U \int_W  \frac{r_n^3}{r_0^3} \left(\prod_{i=1}^{n}\sin \alpha_i \right)
J(a_1,\dots,a_l, b_1, \dots, b_m) \frac{|\xi_n^1|}{r_n} f^{+}(\xi_n)   \{d\alpha_i d\phi_i\}_{i=1}^n.
\end {align} 

We will often refer to this representation in what follows. When we attempt to bound the above expression there will be a factor of $(2\pi)^n$ from the integration in the $\phi_i$ variables and there will be also a factor coming from an appropriate estimate of the $J$-contribution. We want to control these factors. To obtain a good estimate we will split the set $W$ in two parts $W=W_1+W_2$. On the set $W_1$ the $J$-term inside the integral will posses geometric decay with $n$, and on the other hand after a change of variables the measure of the set of integration corresponding to $W_2$ will decay with $n$ sufficiently fast even when $r_0$ is small. We need to do this in order to control $f_n^-$ in such a way that their sum remains an integrable function; if we did not require this, a simpler estimate would still show factorial convergence, but only pointwise in $r_0$. Our extra effort is justified since $f^{-}$ is a probability distribution and it is most natural in addition to pointwise estimates to establish the control of the decay of $f^{-}_n$ in the weighted $L^{1}(\xi)$ norm $\int_{\xi^1<0}|\xi^1| f(\xi) d\xi$ representing the mass flux through the boundary of returning particles. We denote from now on this norm by $L^{1}_w$.

Pick $0< \beta < 1$ and define $k =[\beta n] + 1$ where the brackets indicate the integer part. Also let $D>1$ be a constant dependent on $\beta$ but otherwise fixed that will be specified later. Let $W_1$ and $W_2$ be defined as follows:
$$W_1: (\alpha_1, \alpha_2, \dots \alpha_n) \in W \text { such that } r_k < c/D,$$
$$W_2 = W- W_1 \text { i.e. } W_2: (\alpha_1, \alpha_2, \dots \alpha_n)  \text { such that } r_k> c/D.$$
Notice our definition makes sense since $r_1, \dots, r_k$ are completely determined by $\alpha_1, \dots, \alpha_k$ and $r_0$ so no condition on the $\phi_i$ variables needs to be imposed.  Now with \eqref{fform} in mind we write 
$$ f_n^-= I_1^n + I_2^n$$
where in $I_1^n$ the region of integration in \eqref{fform} with respect to the $\alpha_i$ is $W_1$, and that for $I_2^n$ is $W_2$. To evaluate the integration over these two regions we employ the following lemma. It does not focus on the contribution of the $J$-part of the integrand which will be bounded at the next stage.

\begin {lemma}  The following evaluations hold (assuming $r_0 < c/2$):
\begin {itemize}
\item [i)]
$$\int_{W_1} \prod_{i=1}^n\sin\alpha_i d\alpha_n\dots d\alpha_1 \leq \int_{W} \prod_{i=1}^n\sin\alpha_i d\alpha_n \dots d\alpha_1 \leq \frac{1}{n!}\frac{r_0}{c}\left(\log\frac{R}{r_0}\right)^n$$
\item [ii)]
$$ \int_{W_2} \prod_{i=1}^n\sin\alpha_i d\alpha_n \dots d\alpha_1 \leq \frac{1}{n!}\frac{r_0}{c}\sum_{j=1}^{k} \left( {\begin{array}{c} n \\ k-j \\ \end{array}} \right) \left(\log\frac{RD}{c}\right)^{n-k+j} \left(\log\frac{c/D}{r_0}\right)^{k-j}.$$
\end {itemize}
where we recall that $k = [\beta n]+1$.
\end {lemma}
\begin {proof} 
\noindent i) Pass to the variables of integration $x_i = \cos \alpha_i$. We obtain 
$$  \int_{W} \prod_{i=1}^n\sin\alpha_i d\alpha_n \dots d\alpha_1 = \int_{\Omega} dx_n\dots dx_1$$
where we have $x_i \in (0,1]$ and $\Omega$ with the restriction from \eqref{mainineq} is described as:
$$\Omega: \frac {r_0}{c} \geq x_1x_2\ldots x_n \geq \frac{r_0}{R}.$$
So we need to compute $Vol(\Omega)$. If we further pass to the variables $y_i = -\log x_i$ we get that the corresponding $y_i$ region $\Phi(\Omega)$ is described as
$$ \Phi(\Omega): 0< \log\frac{c}{r_0} \leq \sum_{i=1}^n y_i \leq \log \frac{R}{r_0} \text { and } y_i \in [0, \infty).$$
Furthermore since $x_i=e^{-y_i}$ the Jacobean of the transformation is
$$ \frac{D(x_1,x_2, \dots, x_n)}{D(y_1,y_2,\dots y_n)} = 
Det 
  \left( \begin{array}{cccc}
e^{-y_1} & 0          & \ldots & 0 \\
0        & e^{-y_2}   & \ldots & 0 \\
        & \ldots  & \ldots &       \\
0         & 0& \ldots & e^{-y_n} \end{array} \right) = e^{-\sum_{i=1}^n y_i}.$$
On our region of integration this Jacobean can be bounded as 
$$  e^{-\sum_{i=1}^n y_i} \leq e^{-\log \frac{c}{r_0}} = \frac {r_0}{c}$$
If $\Delta^n_x$ denotes the $n$-dimensional simplex with side $x$ we have that
$$Vol (\Delta^n_x) = \frac {x^n}{n!}$$
so since $\Phi(\Omega)$ is the difference of two such simplices we conclude that
$$ Vol(\Omega) \leq \frac{r_0}{c } [Vol(\Delta^n_{\log \frac{R}{r_0}}) -Vol(\Delta^n_{\log \frac{c}{r_0}})] =  \frac{1}{n!}\frac{r_0}{c}\left(\log^n\frac{R}{r_0}-\log^n\frac{c}{r_0}\right) < \frac{1}{n!}\frac{r_0}{c}\left(\log\frac{R}{r_0}\right)^n $$
as desired.

\noindent ii) Let us perform the same changes of variables as in the previous part of the lemma. Then if $\Omega_2$ is the region in the $x_i$ corresponding to $W_2$ and if $\Phi(\Omega_2)$ is the corresponding region in the $y_i$ from the condition specifying $W_2$ we see that $\Phi(\Omega_2)$ this time is described by
$$
\Phi(\Omega_2): \left\{ \begin{array} {ll}  0< \log\frac{c}{r_0}  \leq \sum_{i=1}^n y_i \leq \log \frac{R}{r_0} \text { and } y_i \in [0, \infty),\\
  0< \log\frac{c/D}{r_0} \leq \sum_{i=1}^k y_i.
\end{array}\right.$$
Notice that the complement in $\Delta^n_{\log \frac{R}{r_0}}$ of the region $\Phi(\Omega_2)$ is described by
$$\Delta^n_{\log\frac{R}{r_0}} - \Phi(\Omega_2) : 
\left\{ \begin{array} {ll} 0< \sum_{i=1}^n y_i \leq \log \frac{R}{r_0} \text { and } y_i \in [0, \infty),\\
 0< \sum_{i=1}^k y_i \leq  \log\frac{c/D}{r_0}
\end{array}\right.$$
Now introduce the unit Jacobean transformation
$$z_j = \sum_{i=1}^j y_i, \text { where } j = 1\dots n$$
so that finally we have
$$ Vol(\Phi(\Omega_2)) = Vol(\Delta^n_{\log\frac{R}{r_0}}) - Vol (complement)= $$
$$ = \frac{1}{n!}\log^{n}\frac{R}{r_0} - 
\int_0^{\log\frac{c/D}{r_0}}\int_{z_1}^{\log\frac{c/D}{r_0}}\ldots\int_{z_{k-1}}^{\log\frac{c/D}{r_0}}\int_{z_k}^{\log\frac{R}{r_0}}\ldots\int_{z_{n-1}}^{\log\frac{R}{r_0}}dz_ndz_{n-1}\dots dz_2dz_1 = $$
$$ = \sum_{j=1}^{k} \frac{1}{(n-k+j)!}\left(\log\frac{R}{r_0}-\log\frac{c/D}{r_0}\right)^{n-k+j}\frac{1}{(k-j)!}\left(\log\frac{c/D}{r_0}\right)^{k-j} =$$
$$ = \frac{1}{n!}\sum_{j=1}^{k} \left( {\begin{array}{c} n \\ k-j \\ \end{array}} \right) \left(\log\frac{RD}{c}\right)^{n-k+j} \left(\log\frac{c/D}{r_0}\right)^{k-j}.$$
As before from the bound on the Jacobean of $\Phi$ we have
$$ Vol (\Omega_2) \leq \frac{r_0}{c} Vol(\Phi(\Omega_2))$$
and so we get the desired expression.
\end {proof}

\begin{subsection} {Bounding the $I_1^n$ term}
On the set $W_1$ we have a strong estimate on the sum in \eqref{fform} as follows. Notice that since $r_k < c/D<c$ we will have that the number of $b$'s in the sums  will be at least $k$, i.e. $ m\geq k $, since all velocities corresponding to $r_0,r_1, \dots, r_k$ must have $\xi^1<0$. Let in fact $b_1, b_2, \ldots, b_{k+1}$ be exactly the ones corresponding to $r_0, r_1, \dots, r_k$ (or in other words to $\xi_0, \xi_1, \dots, \xi_k$). Observe that the requirement $r_k<c/D$ also ensures that for $j \leq k+1$
$$b_j \leq \frac{\pi (c/D)}{c-c/D} = \frac{\pi}{D-1}.$$
 Furthermore notice that we always have $a_i >\pi$, since if $\xi_s$ corresponds to $a_i$ (and so $r_s>c$, $\xi^1_s>0$) then
$$a_i =\frac{\pi|\xi_s-{\bf c}|}{|\xi_s^1|} > \pi\frac{r_s}{r_s -c}>\pi.$$

Now recall that that $I_1^n$ is the same as the expression \eqref{fform} with $W$ replaced by $W_1$. Apply part ii) of Lemma 2.1. with $B = \frac{\pi}{D-1}$ and $A= \pi$ to the $J$ contribution (and bound $b_1 = \pi r_0/|\xi_0^1|$ by itself rather than $B$ which is certainly permitted in the proof of the Lemma). Then after performing the $\phi_i$ integrations we have:
$$I_1^n < \frac{(2\pi)^n}{2\pi^{n+2}} (1+1/\pi)^n \pi^{k}(D-1)^{-k} \frac{\pi r_0}{|\xi_0^1|}\frac{R^3}{r_0^3}\int_{W_1}  \prod_{i=1}^{n}\sin \alpha_i  \frac{|\xi_n^1|}{r_n} f^{+}(\xi_n)  d\alpha_n \dots d\alpha_1.$$
Now we apply part i) of Lemma 2.2., recall that $k = [\beta n]+1$ and bound $\frac{|\xi_n^1|}{r_n}<1$, $f^+_{\xi_n}<M$, to get 
$$I_1^n < M\frac{2^n(1+1/\pi)^n((D-1)/\pi)^{-\beta n}}{2\pi}\frac{R^3}{r_0^3}\frac{r_0}{c}\frac{r_0}{|\xi_0^1|} \frac{1}{n!}\left(\log\frac{R}{r_0}\right)^n.$$
At this point note $2(1+1/\pi)<2.7$ and given the choice of $\beta$ choose $D$ so that 
\begin{equation} \label{dchoice}
\left((D-1)/\pi\right)^{\beta} >2(2.7)=5.4
\end {equation}
which will give finally
$$I_1^n < \frac {M}{2\pi}\frac{R^3}{cr_0|\xi_0^1|}\frac{1}{n!}\left(\frac{1}{2}\log\frac{R}{r_0}\right)^n.$$
This clearly sums to an integrable function near $\xi_0=\bf{c}$ (i.e. near $r_0=0$), in fact we get
$$ \label{ione} I_1 = \sum_{i=1}^\infty I_1^n< \frac{M}{2\pi}\frac{R^{3.5}}{cr_0^{1.5}|\xi_0^1|}\log\frac{R}{r_0}.$$
For the pointwise estimate away from $\xi_0=\bf{c}$, say for $r_0>c/2$, see the remark at the end of this section which applies both to the terms $I_1^n$ and $I_2^n$ and shows that there is no singularity near $\xi_0^1=0$.

The convergence of the $I_1$ term is also fast in the $L^1_w$ norm. If $\mathcal{R}_1^n$ is the remainder term use
\begin{equation} \label{rembound} 
\sum_{i=n}^{\infty} \frac{x^i}{i!}  \leq \frac{x^n}{n!}e^x 
\end{equation}
to see that 
$$\int_{r_0 < R}|\xi_0^1| \mathcal{R}_1^n(\xi_0) d\xi_0 \leq  \frac {MR^3}{2\pi}\frac{1}{2^n n!} \int_0^{R} \frac {|\xi_0^1|}{cr_0|\xi_0^1|} \left(\log\frac{R}{r_0}\right)^n e^{\frac{1}{2}\log\frac{R}{r_0}} 4\pi r_0^2 dr_0.$$
Switching to integration with respect to $u=\log\frac{R}{r_0}$ so $dr_0=-r_0du$ with $r_0=Re^{-u}$ we have:
\begin{equation}\label{firstbound}
\int_{r_0 < R} |\xi^1_0| \mathcal{R}_1^n d\xi_0 \leq  \frac {2MR^5}{c}\frac{1}{2^n n!} \int_{0}^{\infty}u^n e^{-\frac{3}{2}u}du 
 \leq \frac {2MR^5}{c}\frac{1}{2^n n!} n!\left( \frac{2}{3}\right )^{n+1} =  \frac {4MR^5}{3^{n+1}c}.
\end{equation}
\end {subsection}
\begin {subsection} {Bounding the $I_2^n$ term}
For this term we cannot obtain too strong a bound on the $J$ contribution in \eqref{fform}, but by Lemma 2.2. ii) we have a desirable bound on the measure of the integration region after a simple change of variables: namely that measure grows only as a fraction of $n$ power of $\log\frac{R}{r_0}$ which will ensure that when we add the terms $I_2^n$ we will not get a nonintegrable singularity at $r_0=0$. 

As in the remark to Lemma 2.1. since $J^{l,m}<J^{1,1}$ bound the $J$ part of the integrand by $\frac{1}{2\pi}b_1$. This remaining numerator term we can choose to be $b_1$ where we arrange the $b_j$ so that $b_1$ corresponds to $\xi_0$ so $J<1/(2\pi)b_1 = \pi r_0/(2\pi|\xi_0^1|)$ and if we note $|\xi_n^1|/r_n<1$ we get
$$I_2^n  < \frac{(2\pi)^n}{2\pi^{n+2}} \frac{\pi r_0}{2\pi|\xi_0^1|} \frac{R^3}{r_0^3}\int_W  \prod_{i=1}^{n}\sin \alpha_i  \frac{|\xi_n^1|}{r_n} f^{+}(\xi_n)  d\alpha_n \dots d\alpha_1 < \frac{2^n}{4\pi^2} \frac{r_0}{|\xi_0^1|} \frac{R^3}{r_0^3}\int_W  \prod_{i=1}^{n}\sin \alpha_i d\alpha_n \dots d\alpha_1.$$ 

We now study the mass flux due to the remainder term $\mathcal{R}^n_2 = \sum_{i=n}^{\infty} I_2^i$, i.e. we consider the remainder in $L^1_w(\xi_0)$. For this using part ii) of Lemma 2.2. we have
$$\int_{|\xi_0|<R}|\xi_0^1|I^n_2(\xi_0)d\xi_0 < \sum_{q=0, p=n-q}^{q= [\beta n]} \frac{1}{p!q!}
\frac{2^nM}{4\pi^2}\int_0^R |\xi_0^1| \frac{R^3}{r_0^3}\frac{r_0}{|\xi_0^1|}\frac{r_0}{c} 
\left(\log\frac{RD}{c}\right)^p \left(\log\frac{c}{Dr_0}\right)^q 4\pi r_0^2 dr_0 =$$
$$ =  \frac{2^nMR^3}{\pi c}\sum_{q=0, p=n-q}^{q= [\beta n]} \frac{1}{p!q!} \left(\log\frac{RD}{c}\right)^p 
\int_0^R r_0 \left(\log\frac{c}{Dr_0}\right)^q dr_0.$$
In the integral we make the substitution $u = \log \frac{c}{Dr_0}$ so that $dr_0 = -r_0 du \text { and } r_0 = \frac {c}{D} e^{-u}$
and we get
$$\int_0^R r_0 \left(\log\frac{c}{Dr_0}\right)^q dr_0 = \frac{c^2}{D^2}\int_{\log\frac{c}{DR}}^{\infty} u^qe^{-2u} du = \frac{c^2}{D^2} \int_{\log\frac{c}{DR}}^0 + \frac{c^2}{D^2}\int_0^{\infty} = \frac{c^2}{D^2}\left[A_1+A_2\right].$$
By noting that in $A_1$ the integrand is maximized at the left end point and estimating the integral by the maximum of the integrand times the lenght of integration and then computing $A_2$ explicitly we get
$$A_1 < \left(\log\frac{DR}{c}\right)^{q+1}\left(\frac{DR}{c}\right)^2 \text {  and  } A_2 = \int_0^{\infty}u^qe^{-2u}du = \frac{q!}{2^{q+1}}.$$
With this we have
$$\int_{|\xi_0|<R}|\xi_0^1|I^n_2(\xi_0)d\xi_0 < \frac{2^nMR^5}{\pi c}\sum_{q=0, p=n-q}^{q= [\beta n]} \frac{1}{p!q!}  \left(\log\frac{DR}{c}\right)^{n+1} + 
\frac{2^nMR^3c}{\pi D^2}\sum_{q=0, p=n-q}^{q= [\beta n]} \frac{1}{p!2^{q+1}}  \left(\log\frac{DR}{c}\right)^{p} . $$
The first sum is part of a binomial expansion and so
$$\sum_{q=0, p=n-q}^{q= [\beta n]} \frac{1}{p!q!} \left(\log\frac{DR}{c}\right)^{n+1}  < \frac{2^n}{n!}\left(\log\frac{DR}{c}\right)^{n+1}.$$
The second sum is bounded by part of the remainder term in the Taylor explansion of $\frac{1}{2^n}e^{2\log \frac{DR}{c}}$ so as in \eqref{rembound} we have
$$\sum_{q=0, p=n-q}^{q= [\beta n]} \frac{2^p}{p!2^{n+1}}  \left(\log\frac{DR}{c}\right)^{p} < \frac{1}{2^{n+1}}\frac{\left(2\log \frac{DR}{c}\right)^{(1-\beta)n+1} }{((1-\beta)n)!} e^{2\log \frac{DR}{c}}= \frac{1}{2^{n+1}}\frac{\left(2\log\frac{DR}{c}\right)^{(1-\beta)n+1} }{((1-\beta)n)!}\frac{D^2R^2}{c^2}.$$ 
with this we see
$$\int_{|\xi_0|<R}|\xi_0^1|I^n_2(\xi_0)d\xi_0 < \frac{2^nMR^5}{\pi c}\frac{2^n}{n!}\left(\log\frac{DR}{c}\right)^{n+1} + \frac{MR^5}{2\pi c}\frac{1}{((1-\beta) n)!}  \left(2 \log\frac{DR}{c}\right)^{(1-\beta)n+1} .$$
With the above bound on $I_2^n$ we interchange the integration and summation in the $L^1_w$ norm of the remainder  $\mathcal{R}^n_2 = \sum_{i=n}^{\infty} I_2^i$ (which is permitted since the partial sums are bounded by an integrable function as pointed out in the singularity description further below). Thus we have
$$\int_{|\xi_0|<R}|\xi_0^1|\mathcal{R}^n_2(\xi_0)d\xi_0  < \frac{MR^5\log\frac{DR}{c}}{\pi c} \left[\sum_{i=n}^{\infty} \frac{4^i}{i!}\left(\log\frac{DR}{c}\right)^i + \sum_{i=n}^{\infty}\frac{1}{((1-\beta) i)!}  \left(2\log\frac{DR}{c}\right)^{(1-\beta)i} \right] .$$
In both sums above use that $r!>r^re^{-r}$. Thus in the second term use that $((1-\beta) i)! >((1-\beta)i  )^{(1-\beta) i}e^{-(1-\beta) i}$ so that this term is bounded by
\begin {equation} \label{secondbound1}
\sum_{i=n}^{\infty} \left(\frac{2 e\log\frac{DR}{c} } {(1-\beta) i }\right)^{(1-\beta)i} < \sum_{i=n}^{\infty} \frac{1}{e^i} = \frac{e^{-n+1}}{e-1}
\end {equation}
provided that $\left(2 e\log\frac{DR}{c} / ((1-\beta) i) \right)^{(1-\beta)} < 1/e$. To ensure this recall that by \eqref{dchoice} we can choose $D = 1+ \pi 5.4^{\frac{1}{\beta}}$ so $\log{(RD/c)}< 2+ 2/\beta + \log(R/c)$ and so (since $i\geq n$) it is sufficient to choose
$$n > \frac{2e^{1+\frac{1}{1-\beta}}}{1-\beta}\left(2+ \frac{2}{\beta}+\log\frac{R}{c}\right).$$
For the first term similarly choosing $n> 4e^2(2+ 2/\beta + \log\frac{R}{c})$ ensures 
\begin{equation} \label{secondbound2}
\sum_{i=n}^{\infty} \frac{4^i}{i!}\left(\log\frac{DR}{c}\right)^i < \frac{e^{-n+1}}{e-1}.
\end{equation}
Now noting that for $n>2$ one has $\frac{4}{3^{n+1}}<e^{-n}$ we combine the $L_w^1$ estimates on the remainders $\mathcal{R}^n_1$ in \eqref{firstbound} and $\mathcal{R}^n_2$ in \eqref{secondbound1} and \eqref{secondbound2} to get the desired result for $\mathcal{R}_n=\mathcal{R}^n_1+\mathcal{R}_2^n$, i.e. for
$$n>  \left ( 2+ \frac{2}{\beta} + \log\frac{R}{c}\right) \max \{4e^2, \frac{2e^{1+\frac{1}{1-\beta}}}{1-\beta} \}$$
one has
$$\int_{|\xi_0|<R}|\xi_0^1|\mathcal{R}_n(\xi_0)d\xi_0 <   \frac{MR^5}{c}\left(2+ \frac{2}{\beta}+\log\frac{R}{c}\right)e^{-n+1}< \frac{MR^5}{c}ne^{-n+1}.$$ 
If we now choose $\beta = 0.25$ we get the statement of the theorem. We remark again that a more careful analysis may improve markedly the final constants. But for us it was most importnat that these constants are explicit and that we have established the convergence mechanism. This completes the estimates of the mass flux of the remainder terms in our theorem.

Next we show the singularity description. Using part ii) of Lemma 2.2. where we bound the terms in the binomial expansion by the highest power in which they appear we have that 
$$I_2^n  < \frac{2^n}{4\pi^2} \frac{r_0}{|\xi_0^1|} \frac{R^3}{r_0^3}\int_W  \prod_{i=1}^{n}\sin \alpha_i d\alpha_n \dots d\alpha_1 <
 \frac{2^n M}{4\pi^2} \frac{R^3}{r_0^3}\frac{r_0}{|\xi_0^1|}\frac{r_0}{c} \frac{2^n}{n!}\left(\log\frac{RD}{c}\right)^n\left|\log\frac{c/D}{r_0}\right|^{\beta n}.$$
With this we get that the total contribution of the $I_2^n$ terms is 
$$I_2 = \sum_{n=1}^{\infty} I_2^n < \frac{M}{4\pi^2}\frac{R^3}{cr_0|\xi_0^1|}\sum_{n=1}^{\infty} \frac{1}{n!}2^{2n}\left(\log\frac{RD}{c}\right)^n\left|\log\frac{c/D}{r_0}\right|^{\beta n} \text { or }$$
\begin {equation} \label {itwo} I_2 <\frac{M}{4\pi^2}\frac{R^3}{cr_0|\xi_0^1|} \log(RD/c)|\log(c/Dr_0)|^{\beta}e^{C_1(D,c,R)|\log\frac{c/D}{r_0}|^{\beta}}, 
\end{equation}
where $C_1(D,c,R) = 4\log\frac{RD}{c}$. This is an integrable function near $r_0=0$ since as $\beta <1$ for all $\delta >0$ we have 
$$e^{C_1|\log\frac{c/D}{r_0}|^{\beta}}< C_{\delta}r_0^{-\delta},$$
 i.e. the term with the exponent grows slower than any negative power of $r_0$. With this in mind, only for the purpose of the singularity description, we pick on the right hand side of \eqref{dchoice} $3\delta^{-1}$ instead of the choice there (so that in the $I_1^n$ term the exponent of $r_0$ is $-(1+\delta)$ instead of $-1.5$). So together with \eqref{itwo} we complete the singularity description in part ii) of our main theorem:
$$f^{-}(\xi_0) = I_1 + I_2 \leq  C_{\delta} (c,R)\frac{M R^3}{c  |\xi_0^1|} \frac{1}{r_0^{1+\delta}},$$
provided $r_0<c/2$ and where the constant $C_\delta (R,c)$ is explicitly derived from the expression \eqref{itwo}. When $r_0>c/2$ we note that the remark at the end of this section further shows that there is in fact no singularity near $\xi_0^1=0$.

With this the proof of our theorem is completed.

\end{subsection}

\begin{remark} \label{minor}
For the pointwise estimate of $f^-(\xi_0)$ we consider the case $r_0>c/2$ separately because in the above we have $b_1$ for a remaining term in the sum and 
$b_1 = \frac{\pi|\xi_0-{\bf c}|}{|\xi^1_0|}$
is unbounded for $|\xi^1_0|$ near 0 (which can happen if $r_0>c/2$). This point is very easy to correct however since now we are not looking to optimize the singularity at $r_0=0$. In fact we will bound the sum part of the integrand exactly as in the previous section, except this time we choose the term remaining in the denominator to be $a_l$ rather than $b_1$ where we agree that $a_l$ is in fact the term corresponding to $\xi_n$, i.e.
$a_l = \frac {\pi|\xi_n-{\bf c}|}{|\xi_n^1|}.$
Notice that then this term cancels out with the factor $\frac{|\xi_n^1|}{r_n}$ in front of $f(\xi_n)$ in the formula. Thus without splitting the set $W$ this time and just using an argument as the one in part i) of Lemma 2.2. (where the lower bound on the sum of the $y_i$ is now just $0$ so we simply bound the Jacobean of $\Phi$ by $1$) we get for $r_0>c/2$ the bound:
$$f_n^{-} \leq \frac{(2\pi)^n}{2\pi^{n+2}}\pi \frac{R^3}{r_0^3}\int_W  \prod_{i=1}^{n}\sin \alpha_i   M  d\alpha_n \dots d\alpha_1 \leq \frac{1}{2\pi} \frac{R^3}{r_0^3} \frac{2^n}{n!} \left(\log\frac{R}{r_0}\right)^n M$$
The sum of these terms is relevant for $r_0>c/2$ and has no singularity at $\xi_0^1=0$ and vanishes at $r_0=R$ which is the boundary of the domain of influence.
\end{remark}

\begin{remark} \label{Gaussian}
One would hope to extend the above proof to the case when the background is a Gaussian distribution as is very natural for problems in kinetic theory. It is expected that the same basic convergence mechanism will exist in this case but technical difficulties will arise due to the fact that $K$ now will not have the strict shrinking property. In other words due to a diffusive component (in magnitude related to the temperature) the sequence $r_i$ from the proof does not have to be increasing -- however note that this sequence will be close to increasing with a very high probability which may lead to a very similar situation. 
\end{remark}

\end{section}

\begin{section} {Appendix --  Carleman's representation for the gain term}
Consider the gain term 
$$ Q^+(f,g)(\xi) = \int_{\mathbb{R}^3} \int_{S^+} f(\xi-((\xi-\xi_*)\cdot n) n)g(\xi_*+((\xi-\xi_*)\cdot n) n)|(\xi-\xi_*)\cdot n| d\xi_*dn
$$
which according to the geometric intuition explained in the introduction we would like to write as
$$ Q^+(f,g)(\xi_0) = \int_{\mathbb{R}^3} g(\xi_*) W_g(\xi_*,\xi_0) \int_{L_{\xi_*\xi_0}} f(\xi)  W_f(\xi_*,\xi_0, \xi) d\xi d\xi_*$$
where $W_g$ and $W_f$ are weight functions to be determined. We will use the symmetries of the $'$-transformation so we observe
$$ \int_{\mathbb{R}^3} \phi(\xi) Q^+(f,g)(\xi) d\xi = \int_{\mathbb{R}^3}\int_{\mathbb{R}^3}f(\xi)g(\xi_*)\int_{S^+}\phi(\xi')|(\xi-\xi_*)\cdot n| dn d\xi_* d\xi.$$
Now the velocities $\xi, \xi', \xi_*, \xi_*'$ form the vertices of a rectangle and with $\xi$ and $\xi_*$ fixed this rectangle is inscribed in a fixed sphere $S(\xi,\xi_*)$ in different ways according to the value of $n$ ($n$ is parallel to $\xi_* -\xi_*'$). If $d\sigma$ is the surface measure on that sphere one easily finds that 
$$|(\xi-\xi_*)\cdot n| dn = \frac{d\sigma}{|\xi-\xi_*|}.$$
Therefore 
$$ \int_{\mathbb{R}^3} \phi(\xi) Q^+(f,g)(\xi) d\xi = \int_{\mathbb{R}^3}\int_{\mathbb{R}^3}f(\xi)g(\xi_*)\int_{S(\xi,\xi_*)}\phi(\xi')\frac{1}{|\xi-\xi_*|} d\sigma d\xi_* d\xi.$$
and we will take $\phi_\epsilon(\xi) = \frac{1}{(2\pi\epsilon)^{3/2}}e^{-|\xi_0-\xi|^2/(2\epsilon)}$ so that we can compute
$$ Q^+(f,g)(\xi_0) = \lim_{\epsilon \rightarrow 0} \int_{\mathbb{R}^3} \phi_\epsilon(\xi) Q^+(f,g)(\xi) d\xi.$$
After denoting $z$ to be the midpoint between $\xi$ and $\xi_*$ (which is the center of the sphere for the $d\sigma$ integral) a straightforward calculation using the above representation shows
$$ \int_{\mathbb{R}^3} \phi_\epsilon(\xi) Q^+(f,g)(\xi) d\xi =$$
$$= \int_{\mathbb{R}^3}\int_{\mathbb{R}^3}f(\xi)g(\xi_*)
\frac{2}{(2\pi\epsilon)^{1/2}|z-\xi_0|}\left[ e^{-\frac{1}{2\epsilon}\left(|z-\xi_0|-\frac{|\xi-\xi_*|}{2}\right)^2}- e^{-\frac{1}{2\epsilon}\left(|z-\xi_0|+\frac{|\xi-\xi_*|}{2}\right)^2}  \right] d\xi_* d\xi.$$
Notice this symmetric expression shows that $ Q^+(f,g) =  Q^+(g,f)$, but of course $ Q(f,g) \neq  Q(g,f)$ due to the nonsymmetric loss term. It is obvious that the second term in the brackets above does not contribute in the limit $\epsilon \rightarrow 0$. Then a careful but straightforward calculation shows that the first term in the limit is equivalent to $\frac{1}{|\xi-\xi_*|}\delta(h)$ where $h$ is the distance from the plane $L_{\xi_*\xi_0}$ and we indeed arrive at the desired formula
$$ Q^+(f,g)(\xi_0) = \int_{\mathbb{R}^3} g(\xi_*) \frac{1}{|\xi_*-\xi_0|} \int_{L_{\xi_*\xi_0}} f(\xi)  d\xi d\xi_*.$$
\end{section}

\end{document}